\documentclass[sn-mathphys-num]{sn-jnl}

\usepackage{graphicx}%
\usepackage{multirow}%
\usepackage{amsmath,amssymb,amsfonts}%
\usepackage{amsthm}%
\usepackage{mathrsfs}%
\usepackage[title]{appendix}%
\usepackage{xcolor}%
\usepackage{textcomp}%
\usepackage{manyfoot}%
\usepackage{booktabs}%
\usepackage{algorithm}%
\usepackage{algorithmicx}%
\usepackage{algpseudocode}%
\usepackage{listings}%
\usepackage{array}

\begin{document}

\title{Community-driven data science practices}
\markright{Community-driven data science practices}
\author{Atilio Barreda II, Carrie Diaz Eaton, Sam Hansen, Joseph E. Hibdon Jr., Lee T. Gordon, Rebekah Greenwald, María José Gutiérrez Paz, Kenan İnce, Claire Kelling, Drew Lewis, Ariana Mendible, Jenny Mercado, Victor Piercey, Bianca Thompson}

\abstract{
Mathematics researchers are becoming more involved with research questions at the interface of data science and social justice. This type of research needs to be grounded in the needs of the community in order to have significant impact. In this paper, we examine two examples of community-research partnerships in data science for social justice co-authored by both community members and mathematical researchers. The first, VECINA, is a place-based community-research partnership focused on environmental justice. VECINA introduces a framework for developing fruitful local collaborations. The second example, SToPA, originates in citizens' request for an analysis of their town's policing data, but focuses on how to scale this work beyond that place-based setting. SToPA's research helps us imagine how we can  continue to actively collaborate with community members even when working to scale projects beyond a single community. In both of these case studies, we examine the harmonies between established principles of power, process, and perspective with our framework for research-community partnerships. We use a duoethnography approach, directly illustrating the experiences of researchers. We also offer a set of reflections on the impact of these research-community partnerships.}


\maketitle

\noindent

\section{Introduction}
\begin{figure}
    \centering
    \includegraphics[width = .7\textwidth]{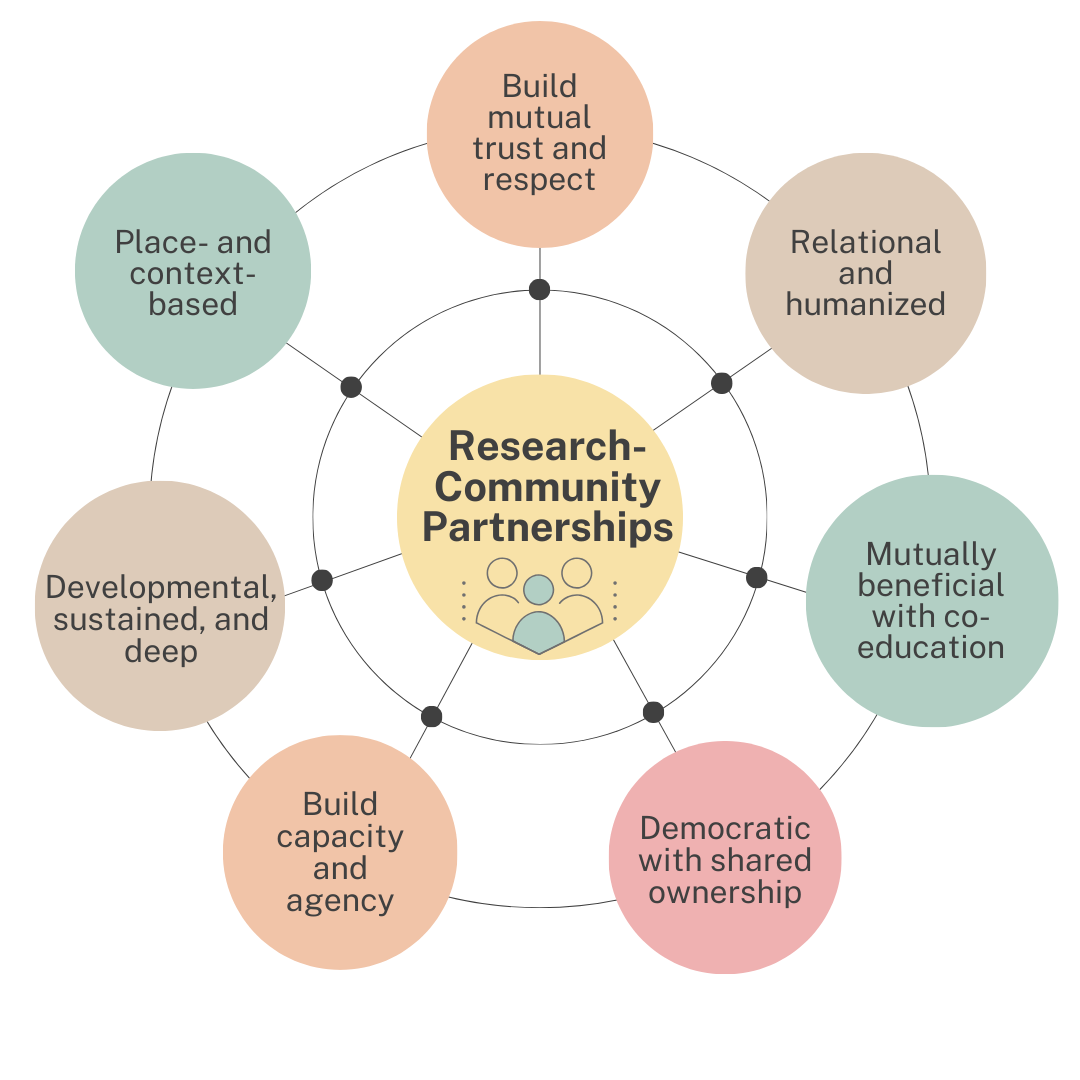}
    \caption{Framework developed during the IMSI research semester for research-community partnerships adapted from the Bonner Framework (\url{www.bonner.org}). We explain how the established principles of power, process, and perspective connect to this framework through the case studies. See Table~\ref{tab:framework} for further elaboration on our framework.}
    \label{fig:framework}
\end{figure}

As conversations of data science and social justice are maturing simultaneously in national discussion, there is more interest in research at the intersection of these fields  \cite{jonesDataScienceSocial2023a}. Research with social justice orientations as well as quantitative methods is common in many fields, such as health and education. In these fields, there is well-developed scholarship on the benefits of and approaches to community-research partnerships \cite{waCommunityBased2010, isCommunityBased2010, racher05}. Our goal is to translate this work into practices for mathematicians who are interested in data science for social justice\footnote{We define data science for social justice as ``data scientific work (broadly construed) that actively
challenges systems of inequity and concretely supports the liberation of oppressed and marginalized communities" as found in \cite{jonesDataScienceSocial2023a} and refer the reader here for a more thorough discussion of the nuances of this definition.}, through the sharing of our own experiences in community-driven research. For example, Figure~\ref{fig:framework} is a framework we created based on our experiences in and lessons learned from these partnerships to give guiding principles to Research-Community Partnerships.

We are a team of mathematicians, community activists, non-profit leaders, statisticians, educators, and computer scientists who worked together on research projects undertaken at the Institute for Mathematical and Statistical Innovation (IMSI) at the University of Chicago in January, February of 2023, and May 2024. The theme of the research collaboration workshop in which we participated was “Interdisciplinary and Critical Data Science Motivated by Social Justice.” Some of us were co-organizers of this program, some of us were applicants to and participants of the research program, and some of us were collaborating on these projects off-site. We hold a variety of social identities and lived experiences - for example, women, men, queer, non-binary, white, Black, Latinx, Indigenous, those with connections to the towns and cities in our research, and those with external perspectives - and we see this heterogeneous mix of professional and social identities as a strength of the work we present here. 

It should be noted that our use of the term “community member” - one who lives in and/or advocates for the community we are partnering with - and “academic researcher” - one who comes to this work enculturated in the academic culture and with a particular disciplinary perspective - does not indicate disjoint identities among our co-authors and research team members. Rather these terms are used in the paper when clarifying the perspective that is invoked when making a claim or using a “we” term.

Members of the communities that are most impacted by issues of social injustice are, in many ways, the most expert on the lived experience and community impact of these issues. Unfortunately, racism and epistemic injustice (injustice related to whose knowledge is considered “valid” or “worth hearing” by those in power \cite{kiddRoutledgeHandbookEpistemic2017}) mean that those who have lived experiences with injustice are often not seen to have the epistemic standing (“enough knowledge or expertise”) to challenge existing power structures in which government and academic narratives are considered authoritative. Since data science and mathematics, especially when done by academics, are often considered “objective” (as in the common saying “Numbers don’t lie”) data science can be used to provide a quantitative perspective on the lived experience of communities most vulnerable to injustice in a way that is more likely to be respected by people in power \cite{symonsEpistemicInjusticeData2022}.

As research-community partners, we focused on issues of social injustice, specifically related to issues of policing, climate change, racial justice, and environmental justice. Each of these projects has a slightly different origin story. The first project - VECINA, was initiated as part of another research semester experience at the Institute for Computational and Experimental Research (ICERM) at Brown University in Providence, RI in 2022. The research-community partnership started at ICERM, and this experience formed the basis of the research undertaken at IMSI in 2023 and has evolved into many VECINA teams (equipitos) working in tandem. The second project, the Small Town Policing Accountability (or SToPA) Lab, has its roots in the work of the Institute for Quantitative Study of Inclusion, Diversity, and Equity (QSIDE Institute). The SToPA Lab was developed in 2021 by a group of data scientists, social justice activists, and statisticians in response to an epidemic of white supremacist practices within the Williamstown, MA police department. The lab was later expanded to include Durham, NC, at Duke University to create transparency and equity in policing; therefore also building capacity for restorative justice. In Spring 2024, part of SToPA's work has led to the formation of a new research group focused on communities of all sizes called the Data Science, Police Accountability, and Community Engagement (DSPACE) Group. This group intends to make a direct and meaningful impact with partner communities as well as to contribute to scholarly research in data science and statistics. The work of both VECINA and SToPA at IMSI was part of a broader effort to integrate data science research into the social justice struggle\footnote{Data science for social justice (DS4SJ) is data scientific work (broadly construed) that actively challenges systems of inequity and concretely supports the liberation of oppressed and marginalized communities.\cite{jonesDataScienceSocial2023a}} to help activists in other jurisdictions.

A recent article in the \textit{Notices of the American Mathematical Society} suggested a framework for approaching data science for social justice research, especially given the increasing interest sparked in 2020 \cite{jonesDataScienceSocial2023a}. The first principle of \textbf{\textit{power}} attends to the need for social justice research to have an impact. As a result, researchers need to consider and identify who has the power to create change. This can include identifying possible collaborative community partners and audiences for the research. What matters to community audiences is crucial in formulating research questions and identifying data needs or knowledge partners. The second principle of \textbf{\textit{process}} emphasizes humbly engaging with knowledge partners and community partners to promote impactful long-term relationships rather than transactional ones. The third is \textbf{\textit{perspective}}. Knowledge partners include those in disciplines that can help inform social justice work, such as those in the social sciences, humanities, and/or fields of application (e.g., community health or environmental studies). However, this final principle emphasizes that the insight of direct, lived experiences is key and often distinct from disciplinary knowledge. 

This framework exists alongside other frameworks for justice-based community-engaged partnerships, each focusing on best practices for a specific audience seeking engagement with community partners. The Anti-Racist Community Engagement framework was formulated for institutions of higher education \cite{santana_anti_racist_2023}, the Design Justice Framework for technologists \cite{designjustice}, and the Kirwin Institute's Principles for Equitable and Inclusive Civic Engagement \cite{Kirwin} and Bonner's Community-Based Research framework \cite{bonner_community_based_2023} for more general community partnerships. red Our framework focuses on the particular intersection between quantitative academics and community partners and is built on the lessons learned in our partnerships. This is in no way a complete list, and as societies, cultures, and communities change and evolve, so will the methodologies. Nonetheless, we believe these frameworks, as well as our own, are important starting points in scientific inquiry.

Many of us came to mathematics and social justice by first examining social justice in our classrooms or in the mathematics community. In these cases, we advocate for justice in academia, so we become both researchers and target community members. But as we look to support social justice outside of the academic community, this means that we are challenged to engage with a broader set of partners that can bring a critical voice to the process and guide the project to its maximal impact. 

In this paper, we show how two projects, VECINA and SToPA, integrated the three principles of power, perspective, and process \cite{jonesDataScienceSocial2023a} into the research process and how these principles are embedded in our framework shown in Figure \ref{fig:framework}.  We also explain the products and relationships that were developed in both teams. We also discuss the impact of these two projects as well as challenges and lessons learned. These case studies are presented as duoethnography, which allows individual perspectives to be heard while also synthesizing those perspectives to create new and emergent insights across projects \cite{saDuoethnography2012}.

\section{Nuevas Voces, WRWC and VECINA}
Co-author members: AB, CDE, JEH, RG, MJGP, DL, JM, VP, BT
\subsection{Perspective}
\textbf{\textit{What is the context? Who in the community is working with and identifying the changes needed? Who in the community is working to make that change happen?}} \textbf{\textit{How did the knowledge of all research team members [including community members] guide the project in various directions?}}

\begin{figure}[h]
\includegraphics[width = \textwidth]{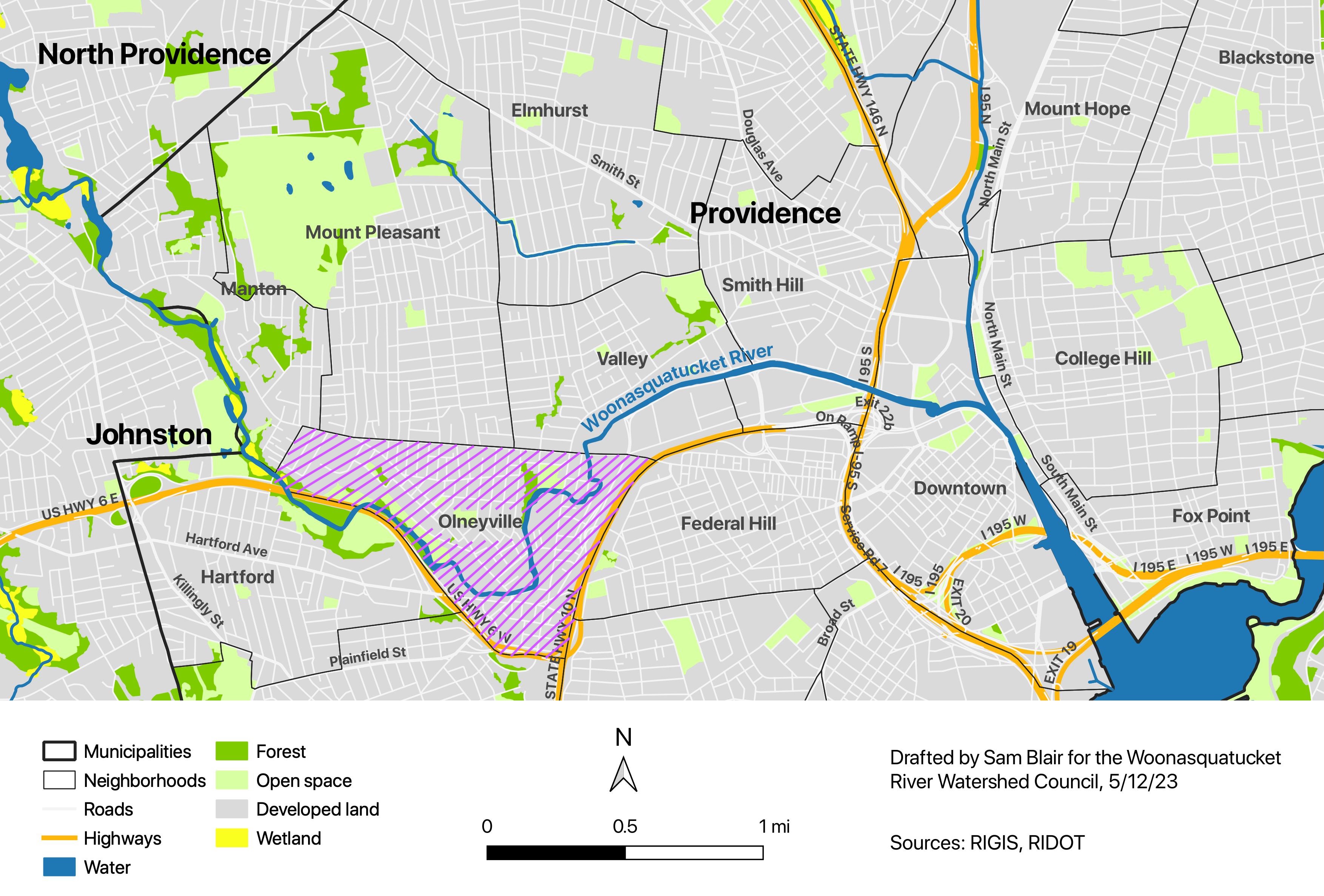}
\caption{The Woonasquatucket - Algonquian name meaning “where the salt water ends” - runs through downtown Providence, RI and into the Narragansett Bay. Brown University (in College Hill) and ICERM (in Downtown, just west of College Hill) are located to the east of the river, while Olneyville (highlighted with purple diagonal lines) is a neighborhood bordering the river to the west of the bay.}
\label{fig:Providence_Map}
\end{figure}

The Woonasquatucket River runs through downtown Providence, Rhode Island as seen in Figure \ref{fig:Providence_Map}. The Woonasquatucket River Watershed Council (WRWC) is a non-profit located in Providence that focuses on restoration of the Woonasquatucket and its watershed area, greenway and greenspace construction and revitalization, and fostering community around the river. In 1999, it was discovered that the Woonasquatucket River was one of the most dioxin-contaminated rivers in the country due to the mills on the river \cite{zaToxicity2002} and mobilizing activism has faced a number of hurdles in its response \cite{zaIssue2004, zaToxicity2002}. The mission of WRWC since its founding in the 1990s is to create positive environmental, social and economic change by revitalizing the Woonasquatucket River, its greenway and its communities, particularly in the downriver neighborhoods most subject to severe storm flooding, heat island effects, and other legacy and climate challenges. Olneyville is one of the downriver neighborhoods which borders the Woonasquatucket (see Figure \ref{fig:Providence_Map}). 

\textbf{Place and context-based} Historically, Olneyville has been a primarily rental neighborhood for many immigrant families, a majority of which are from Central America and the Caribbean. In addition to issues of pollution, the Olneyville neighborhood has experienced significant flooding due to weather extremes from climate change and poor infrastructure \cite{jerzyk09}. In 2020, WRWC created a new program, “Nuevas Voces” or “New Voices at the Water Table” to support emerging community leaders in Olneyville - investing in key community partnerships important to its success. Nuevas Voces was designed for and with the Olneyville community. After the first year, two Nuevas Voces Fellows, co-authors JM and MJGP, began as co-facilitators in the next iteration of the program.
\begin{quote}
\textbf{JM}: I have lived much of my life in Olneyville and all the time I was focused on the day to day. In 2021, I got involved in the Nuevas Voces program where I learned a lot, and it opened my eyes. I learned how to take action on what was happening around me. One of the great things I learned was to prepare for an emergency, but also how to use my voice to help improve my neighborhood and my community. 

For me, helping is something that comes from my heart, and this has led me to want to learn more, and thus to be able to help and give the necessary tools to those who need them most. That's how I decided to be part of facilitating the Nuevas Voces classes because I want to see more people ready to use their voices and be able to help improve our community and make a better future for our children.
\end{quote}

\begin{quote}
\textbf{MJGP}: El amor me trajo a E.E.U.U., al país reconocido por los inmigrantes como el país de las oportunidades. No entendía esta frase... hasta que la salud me jugó una mala pasada y finalmente, aquí estoy, sin hablar muy bien el inglés y aunque felizmente casada después de recuperar mi salud, tuve que comprender que mi vida acababa de dar un giro de 360 grados. Sentí que debía buscar oportunidades para abrirme mis propios caminos y conocí Nuevas Voces, el programa que me brindó conocimientos que yo necesitaba y me sentí bienvenida (no discriminada), empoderada y apoyándome en personas capaces de escucharme y ayudarme a adaptarme al sistema del país, empecé a desenvolverme como líder de mi comunidad para elevar mi voz y la de aquéllos que no lo hacen y que viven con temor por diferentes razones: no saber el idioma inglés, no haber estudiado, no entender el sistema, no ser ciudadano americano y sencillamente, por venir de ciudades o campos humildes de latinoamérica con una realidad muy diferente.

Así mismo, el programa de Nuevas Voces, me hizo sentir como en mi hogar, porque en mi ciudad Santa Cruz de la Sierra, región amazónica de Bolivia, vivimos rodeados de un río y la problemática de las inundaciones es algo similar, sin embargo nunca tuve un curso sobre resiliencia climática y Nuevas Voces me ha permitido despertar a ese importante hecho activándome para ser capaz de manejar situaciones de emergencias, ayudándome a mí, a mi familia y a mis vecinos.

\textbf{MJGP} (translation by JM): Love brought me to the United States, the country recognized by immigrants as ``the country of opportunities". I didn't understand this phrase... until my health played a trick on me and finally, here I am, not speaking English very well and although happily married after recovering my health, I had to understand that my life had just taken a 360 degree turn. I felt that I had to look for opportunities to open my own paths and I met Nuevas Voces, the program that gave me the knowledge that I needed and I felt welcome (not discriminated against), empowered and supported by people capable of listening to me and helping me adapt to the country's system. I began to develop as a leader in my community to raise my voice and that of those who do not and who live in fear for different reasons: not knowing the English language, not having studied, not understanding the system, not being an American citizen and simply, for coming from humble cities or fields of Latin America with a very different reality.

Likewise, the Nuevas Voces program made me feel at home, because in my city Santa Cruz de la Sierra, in the Amazon region of Bolivia, we live surrounded by a river and the problem of flooding is somewhat similar, however, never I had a course on climate resilience and Nuevas Voces has allowed me to wake up to that important fact by activating me to be able to handle emergency situations, helping me, my family and my neighbors.
\end{quote}

To understand why it is important to engage the local community in environmental justice advocacy, we need to understand the ways that Providence is spatially segregated. For example, contrasting Olneyville and the Brown University neighborhood of College Hill (see Figure \ref{fig:Providence_Map}) illustrates how the “Ivory Tower” effect results in segregation based on education, race, and ethnicity \cite{laBuilding2017}. In College Hill 80.2\% of residents have a college degree and 71.4\% are non-Hispanic White, whereas 16.8\% of Olneyville residents have a college degree and only 16.5\% are non-Hispanic White \cite{Demographic}. This is accompanied by economic segregation, especially among the wealthiest residents. In College Hill, the top 5\% income households make an average of \$619K per year, which is the highest of all Providence neighborhoods. In contrast, the top 5\% of Olneyville households make less than \$124K, which is the least of all Providence neighborhoods \cite{Demographic}. College Hill, as indicated by the name, sits on a hill and encloses one repeated flood site, in contrast with Olneyville which includes land in the 100- and 500-year FEMA floodplain and multiple repeated flooding sites \cite{Providence2021}. This comparison between College Hill and Olneyville shows that climate change is also a social justice issue. 

VECINA, Visualizing Environmental and Community Information for Neighborhood Advocacy, was the result of a joint project between researchers at ICERM and IMSI and research partners and community leaders at WRWC and Nuevas Voces. The statistics shared above reveal why research about social justice at ICERM at Brown University is not the same as local community-engaged social justice research. Not only are the identities very different between the university's neighborhood and nearby local community members, the lived experiences brought to the research process also often differ severely. The academic researchers would have to leave the ivory tower to seek relationships which foster social justice work. Throughout this section we account the collaborative process through our team’s voices and summarize the key tenants (Table \ref{tab:framework}) which we hope will guide others’ work. We also connect this collaborative process and framework to the principles of power, process, and perspective, through our local and placed-based justice work.

\begin{table}[h]
    \centering
    \begin{tabular}{
    | >{\raggedright\arraybackslash}m{2cm}
    | >{\raggedright\arraybackslash}m{8cm} 
    |}
    \hline
     \textbf{Collaboration feature} & \textbf{Description}\\\hline
        Place and context-based & Understand where the work is grounded, the community involved, and the broader history as they should inform the project, the collaboration, and the customized approach.\\ \hline
        Build mutual trust and respect & Develop a community of care in which you respect each other's personhood, knowledge, and expertise - then sow the seed of deep commitment by having open communication and counting on each other to follow through. \\ \hline
        Relational and humanized & This work is about people, relationships across researchers in academia and community as well as relationships between communities within the broader sociopolitical context. \\ \hline
        Mutually beneficial with coeducation & Each person and partner organization can contribute their diverse sets of experience and expertise and has identified the ways in which they will benefit and/or grow from the broader collaborative effort.\\ \hline
        Democratic with shared ownership & There is a transparent and clear process for making group decisions about the project and as a group and a shared understanding of project ownership - both should be guided by community needs. \\ \hline
       Build capacity and agency  & Sustainable change requires that we build our knowledge and capacity to make change in all of us as individuals and organizations, and as such we should help each other grow into this potential. \\ \hline
       Developmental, sustained, and deep  &  This journey is ongoing, and each member and organization should be transparent about its commitment to each other, the project, and the broader vision. Make space with open communication for this to change and evolve over time for sustainability.\\ \hline
       \hline
    \end{tabular}
    \caption{Definitions for the framework presented in Figure~\ref{fig:framework}.}
    \label{tab:framework}
\end{table}

\begin{quote}
\textbf{CDE}: I was an organizer of the semester research program at ICERM. I applied to ICERM as a host the social justice for data science program because of the location of the Institute’s home university at Brown in Providence, RI. I grew up only a few miles from ICERM’s location, and my Latinx community was based out of a Spanish-speaking church near Olneyville. It was this church community that was at the heart of the community support for Latinx immigrants and this community helped build my identity as a Latin-American in a predominantly white state. I wanted to initiate a project that was giving back to this while I was in residence at the institute. 

It had been years since I lived in RI; and it was difficult in particular to make connections, because we were just moving out of the COVID era. Many non-profits were struggling and did not have the bandwidth to take on discussions that were nebulous. But we leveraged individual connections, colleagues and friends of friends who could connect us to potential non-profit partners. These connections are often the most fruitful because they benefit from the previous relationship work of others and there is a brokerage of these relationships based on trust and reciprocity. In RI in particular, however, there is a culture of “I know a guy who knows a guy” - this often spoken phrase is indicative of a small world brokerage network which strongly operates in this context.

My sister still lives in RI and works in the environmental field. She was the one who introduced me to WRWC. She was raving about their new Nuevas Voces project and the work they had undertaken in the Olneyville community restoring the watershed with the help of small grant projects and advocating for their community. She set up a lunch meeting with the director of the WRWC and our team to get to know each other. 

It was clear that Nuevas Voces was a program that was making a difference in the community. But neither of us knew what it might mean for mathematicians to support that work. However, through a Data Science for Social Good training program \cite{dhData2021}, I learned to ask the question of partners: “What case do you have to make over and over again, that having the data to back you up would be helpful?” The phrasing of this question highlights that they are the change makers and we as mathematicians are there to support them with a set of skills we can bring to the table. The meeting with WRWC led to some interesting ideas, but it was clear that the next step was to engage with the Nuevas Voces leaders who really understood the boots on the ground needs of the community. 
\end{quote}

\textbf{Mutually-beneficial with co-education} 

We collectively co-educated each other in various ways throughout this process. 
Each participant brought unique skill sets to the table - Spanish-English translation skills, R and data cleaning/ visualization skills, knowledge of Providence community needs, knowledge of how nonprofits can leverage VECINA, knowledge of databases and back-end development, and more. At IMSI, we also had access to the skills and connections of those at the SToPA Lab and the Center for Spatial Data Science at UChicago, which housed an open-source spatial visualization tool, Chives (https://chichives.com/). 

\begin{quote}
\textbf{CDE}: Before arriving at IMSI, we already knew we wanted to create VECINA to combine stories with geospatial information. However, we initially thought the emphasis would be on environmental data, given the role of WRWC as interested in climate change and watershed protection. However, JM and MJGP suggested developing a community resource list first. We also were prepared to have the site in Spanish, but both advocated to also have an audio button because many immigrants don’t have literacy skills - not because of general literacy, but that many immigrants come with language skills in Indigenous languages, not the colonizer language of Spanish. While MJGP was still improving her English, she was fluent in many other languages. Also, JM talked about her current frustrations with how to fill out school choice preference forms. She really helped us understand a multitude of other issues that were significantly impacting the community and this led to conversations resulting in the addition of another layer to VECINA that would help parents with choosing schools.
\end{quote}

\begin{quote}
\textbf{JM}: I believe that education is the key to a better future full of knowledge. The Providence educational system is not designed to help parents or students. We're almost letting them fail. My experience in helping to choose the next school for my nephew was frustrating because we did not have any information on what was the best school for him, and this surely happens with thousands of families in Olneyville. Not having access to information makes it more difficult to choose the school according to the needs of the student and the family. During my time in Chicago I had the opportunity to learn from the team that my nephew didn't have lots of choices, that none of the schools offer him access to bilingual classes, or the academic goals that we were looking for him to pursue. VECINA will be able to help give us this information and thus a better understanding of which is the best option - to know if the school offers bilingual education, the distance, if they are reaching the goals of teaching and educating the students. For me, education is also about justice because if we fail our children, we are condemning them or not giving them a future.
\end{quote}

\subsection{Process} 
\textit{\textbf{How was the process of collaboration designed to enact the collaborative framework?}}

\textbf{Relational and humanized} Dr. Sunshine Brosi, in a plenary lecture at the BioQUEST Summer Institute in 2018, shared that when she starts teaching students to do community-engaged work, she teaches them first how to make a pie. The pie was an offering of time, effort, and learning that meant you were serious about putting in the time to relationship build and wanted to signal that you are treating them like a neighbor. It was a moment in which many in the audience reflected on how and why we do not acknowledge this step. It reminds us that academics have to approach communities with respect and humility. We have to be willing to give something of ourselves that is valuable - time and money - to communities and community members. 

The academic researchers of ICERM left our ivory tower, brought food, and visited WRWC headquarters. We set aside ample time to learn about each other and communicate each other's goals and the needs of WRWC and Nuevas Voces. The meeting was conducted in Spanish, with translation into English as needed. We spent time walking along the river with our guides from Nuevas Voces, hearing stories and seeing with our own eyes where flooding had happened and where restoration was happening. It was at this meeting that we identified possible avenues of mutual interest. We continued to meet monthly as we formed VECINA and found ways to support the work through IMSI as well as a grant from the [redacted] Center, the Community Engagement office at CDE's home college, to support stipends for Nuevas Voces leaders. Once we knew that we were going to take on the work at IMSI, co-authors JM and MJGP worked to develop a wishlist that could kick start the conversation at IMSI.

\begin{quote}
\textbf{AB}: In developing software solutions for community-centric projects, the principles of adaptability and experimentation prove to be pivotal. This echoes familiarly within the realm of academic mathematics, where these same principles enable us to delve deeper into unknown territories, uncover new theories, and build upon existing knowledge. The approach we take transcends mere technology and reflects the inherent diversity and richness of the community we serve.

At the heart of this process is the recognition of community organizers as product managers and stakeholders—an acknowledgment that invites them to be active contributors in the development process. It is a shift from a transactional to a relational approach, fostering a partnership that is rooted in trust, mutual respect, and a shared vision for positive change. 

In this context, the software development process becomes a medium for storytelling—a platform where individuals like JM can share their experiences and aspirations, and use these narratives to shape the design and functionality of the software. It is a process where every line of code, every user interface, and every feature is imbued with the values, experiences, and aspirations of the community.
\end{quote}

\textbf{Building a team of trust and mutual respect} From the WRWC and other activist group perspectives, Brown is a university with a significant amount of resources. However, many city improvement and environmental impact projects are designed to financially benefit Brown and contribute to academic knowledge, leaving Providence’s frontline communities out of investment, power sharing and problem solving. These are very different outcomes from those of Nuevas Voces, which is being designed to support the community. This is why the idea of “with not for” - or “nothing about us without us” as popularized by the international disability justice community - is so important in this work \cite{jonesDataScienceSocial2023a}. This is why relationship building as a team together is so important, and why sometimes first steps must include repairing broken trust \cite{waCommunityBased2010}.

\begin{quote}
\textbf{JEH}: At IMSI 2023 we were given the opportunity to work directly with people of Nuevas Voces and WRWC to develop VECINA. Nuevas Voces leaders generously took time to come work with us. Before getting together Nuevas Voces engaged with the community on what exactly was needed for this online app that would help meet their needs. 

There was no guarantee we could meet all the needs that were requested but through conversations with Nuevas Voces we were able to narrow down the most important items based on the overall conversation with the community. The development of who can and was willing to work on what items was a collaborative conversation between Nuevas Voces and the rest of the individuals working on the project. At this point it becomes a little redundant and not fair to state that Nuevas Voces, WRWC, and everyone else as being separate but as a team. 

At the end of each day we shared the outcomes of the work and made sure if anyone wanted to know more information that we would teach each other what we had found.  Development of continued communication was key in getting a good product, but it was also important since what those not part of Nuevas Voces deem important might not be what the community of Providence really wants.
\end{quote}

\begin{quote}
\textbf{JM}: Being part of this research taught me how important it is to use data and mathematics and to be able to use these to give more truth to the problems facing our community. Only with specific data and numbers can we accurately demonstrate the changes we need to make, the things we want to change or achieve in the long term.

The main takeaway for me is that together we can achieve more things for the community and reach those who feel that their voices can't be heard because of language barrier or the lack of knowledge of the differences in resources. 
\end{quote}

\begin{quote}
\textbf{RG}: For WRWC, the launch of Nuevas Voces in 2020 was intended to repair community trust. Our hope is to become truly community-led with watershed residents – particularly those from the downriver neighborhoods like Olneyville where community members have been historically marginalized – in positions of leadership across our sector, our organization, our collaborations with partners, and as the lead decision makers in problem solving where we work. As with VECINA, “with not for” is a pillar of Nuevas Voces, Campeones de Combate Climático and a newly developing resilience hub also being designed by Olneyville residents to meet their communities’ needs.

Personally, I feel grateful to learn through the VECINA process to become an even stauncher and possibly more effective ally than I already thought I was. I feel both humbled and privileged to work with leaders like JM and MJGP, and other coworkers, participants and graduates from the neighborhoods.
\end{quote}

\subsection{Power} 
\textbf{\textit{What aspects of the projects were designed to provide communities with the power to take on this work? What additional capacity is enabled by this work?} }

\textbf{Democratic with shared ownership} As a team, we research possible platforms, data sources, community needs, and more and figure out together what would serve Nuevas Voces the best, including thinking about ownership and maintenance early. Nuevas Voces had already compiled a list of resources that they wanted to publicize. A prior intern had developed an ArcGIS StoryMap that we re-discovered on the internet after a Google search. The StoryMap was limited in what it did other than plot names and locations. No one had remembered it, and no one had access to modifying it. This is why we prioritized using open-source software and drawing from databases owned by the organizational WRWC Google Drive. We wanted to make sure WRWC researchers were the owners of data used by VECINA and they had the ability to manage the data themselves.

When we wanted to name our version of a story map together, we came up with the name “Fuertes de la comunidad,” which translates to “Community Strengths.” This language reflects a positive view of community resources. Throughout the process at IMSI we co-designed this project together, prioritizing Nuevas Voces needs.

In our work, we center the open ethos of collaborative and democratizing philosophies that underlie open science, open-source technologies, and open education. Open-source software allows for free public use and contribution. Within their product communities, users can also become contributors. Therefore, choosing open-source software not only allows for accessibility, but also allows us to act both as creators and consumers, taking insights from the community's collective wisdom, and offering our innovations in return -- a reflection of the cooperative spirit within our research group. 

\textbf{Build capacity and agency} Academic researchers on the team took on projects that developed data science skills and many thought about how to bring this type of experience back to their classrooms and research labs. Also, as we were coworking, Nuevas Voces and WRWC researchers were adding folders and file structures, so they were familiar with where the folders were and how to edit to add new places or modify translations for the StoryMap (and now webpage).

As mathematicians, engineers, data scientists, and community organizers we traverse a fascinating realm of collaboration and situated knowledge\footnote{Dona Haraway's ``Situated Knowledges" presents a foundational perspective for understanding the relevance to data science practices. Haraway argues for an understanding of knowledge as always partial and located within specific perspectives and contexts. This resonates with social justice interventions in data science, as data is never neutral but inherently influenced by the conditions of its collection, interpretation, and application. In practice, this implies that data scientists must account for these factors, interrogating not only the data itself but also its context and the perspectives shaping its use.} In this technological assembly, we find parallels with our interdisciplinary research - selectively choosing, adapting, and integrating tools. Rather than adopting a uniform approach, we utilize a rhetorical and adaptive approach to technology selection. We choose to employ elements in response to community infrastructure, community needs, affordability and accessibility, and future project growth.

\begin{quote}
\textbf{RG}: A world of opportunity has opened up to democratize data since we started this project. Coming into the 3rd year of Nuevas Voces, as a community we are thinking about data differently, I would say more in a justice framework, and we are building and using data differently in the program. As the site comes online, we know that we still have a lot to learn related to populating and running the site. Luckily, our research partners have been very generous in helping us with our efforts to secure additional resources so we can train, plan for widespread community engagement in visual, oral and written storytelling, and form a work group who will carry the project forward, maintaining its integrity as a community space. 
\end{quote}

\subsection{Products and Relationships}
\begin{quotation}
“Move at the speed of trust”
-adrienne maree brown (pg. 16, 2021) \cite{adHolding2021}
\end{quotation}

\textbf{Developmental, sustained, and deep} We pair “products” and “relationships” for three reasons. First, products produced should be contextualized by how they serve the community and enhance the community-research partnership. Second, relationships are an extremely valuable “product”, though we use quotes to intentionally reject the capitalistic urge to quantify the value of a relationship. Third, we suggest that the most important outcome is a genuine relationship, because innovation and products will derive from a collaborative relationship, but the reverse is not true.

After the conversations in summer 2022 at Providence, a “product” was an agreement to collaborate based on a shared understanding between the mathematics research team and WRWC about what we were working toward and how it would benefit both parties. While VECINA was not fully specified yet, we had a vision that we would work together to create some kind of community information platform around environmental justice. Throughout the fall, we negotiated what that collaboration would look like, who would participate, and the roles, and we coordinated our time at IMSI. 

During the winter at IMSI, we collaborated to create specifications for VECINA and priorities. We also launched a number of smaller working groups around each piece of work to address the highest priorities. For example, one major challenge became selecting an appropriate software for the environmental justice information platform. Incredible database and visualization tools, like EJScreen (\url{https://screening-tools.com/epa-ejscreen}), did not have the UX that was useful. An unofficial copy of EJScreen illustrates this tool now, as the original has been deactivated on the government website \cite{EPAEJScreen}. Chives was open source and had a great format for visualizing spatial information, but it did not have functionality to add point data. Leaflet makes great StoryMaps with point data, is open source, and is able to draw from Google Spreadsheets. No programs were set up with Spanish and English languages and audio. Because of these limitations, we determined that multiple platforms should be integrated to achieve the optimal solution for this use-case. Therefore, the VECINA product at this stage was actually a series of proof-of-concept or beta test individual products, to be later integrated. 

While we were at IMSI, JM shared that she was advising a new immigrant family on choosing a middle school. Providence is a hybrid school choice, meaning that students can go to a neighborhood school within one mile or apply to attend a different school in the city. As the prospective middle schooler did not speak much English, JM shared that there is limited information about which schools offer special English Language Learning (ELL) programs beyond what is legally required for accommodation. This led team member VP to start investigating test outcomes at all schools between ELL and non-ELL students, and ultimately creating a dashboard for families to learn about local school options. 

\begin{quote}
    \textbf{VP}: JM's question is not only important for the family she is advising, but for all of the parents in her community. This is especially true for new immigrants who are still acclimating themselves and may not have access to``word-of-mouth" information about schools. In our initial research, we found that the school choice process in Providence was very opaque. Parents often consent to send their children to a school within one mile of their home. The work we have done now helps parents make more informed decisions about where to send their students.

    These data can do more than just inform school choice preferences. They can form a basis for the community to make policy recommendations about where to place educational resources. For example, a community may want to persuade policymakers to invest more in English Language Learning programs in schools closer to their neighborhood, or in schools with higher test scores.
\end{quote}

\begin{figure}
    \centering
    \includegraphics[width = \textwidth]{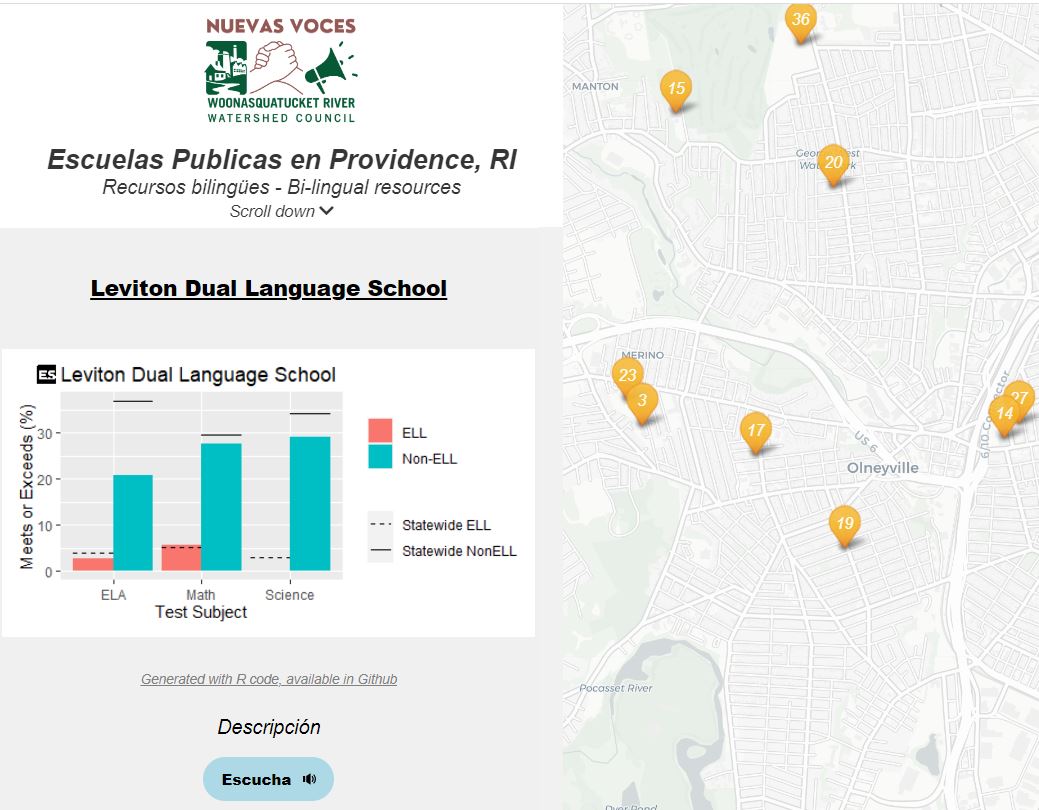}
    \caption{Screenshot of Leaflet StoryMap showing testing outcome disparities between ELL and non-ELL students. The ES symbol in the upper left by the school name indicates that this school has specialized programming for ELL students.}
    \label{fig:pvd-edu}
\end{figure}

The resulting Leaflet StoryMap can be found here \url{https://math-data-justice-collaborative.github.io/edu-pvd/}
and is linked to in the VECINA website hosted by WRWC (\url{https://vecina.wrwc.org}). A screenshot is provided in Figure \ref{fig:pvd-edu}. After reviewing the analysis of all schools, we found tragically depressing results. In addition to the overall low test-passing rate in Providence public schools, there was a huge disparity between ELL and non-ELL students. While our tool helps parents make choices, there seems to be no good choice. In 2019 prior to our work, exposure of such discrepancies led to a state takeover of the Providence schools \cite{RItakeover}, but this takeover has faced many challenges from legislation limiting authority to COVID-19.

Through this process, the “shared understanding and willingness to collaborate” matured into a “relationship based on mutual trust and respect.” By letting the community members lead the goals and use of our products, resulting tools become more useful to those who need them. By letting the mathematicians lead what is possible with the data and technical skills, they can help drive insight. This dynamic reciprocity looks like co-leadership, but builds us as a community together. In the second summer at ICERM, team members continue to engage on site, with an emphasis on integrating issues of health disparities into existing tools. As seen through the eyes of Nuevas Voces, environmental justice is not separate from educational justice, health justice, and social justice. 

\section{SToPA}	
Co-author members: AM, CK, Kİ, LTG, SH 

The Small Town Policing Accountability (SToPA) Lab is a research group within the Institute for Quantitative Study of Inclusion, Diversity, and Equity (QSIDE Institute) \cite{instituteforthequantitativestudyofinclusiondiversityandequityQSIDEInstitute}. The primary goal of the SToPA Lab is to facilitate the analysis of police data in partnership with local community groups in smaller towns and cities. Community groups have long protested against police misconduct and conducted analyses to document racial bias in policing. Since the development of policing in the U.S. from Southern slave patrols 
\cite{kalabhattarHistoryPolicingUS2021}, many community members and activists have been concerned with the racial dynamics of policing throughout the country. 

Nationally, racial bias in police use of force gained renewed attention  in the national conversation in 2014 with the killings of Trayvon Martin, Michael Brown, and Eric Garner, all unarmed Black men killed by police. One year later, the Washington Post found that unarmed Black men were seven times more likely than whites to die by police gunfire \cite{somashekhar2015black}.
More recently, the protests beginning in May 2020 that resulted from the police murder of George Floyd in Minneapolis, MN were the largest in U.S. history, involving between 15 and 26 million people nationwide \citep{buchanan2020black}.  The global circulation of video documenting George Floyd's murder led to protests that spread to over 2,000 cities \citep{buchanan2020black}. 

Much of the data science work on police injustice is produced by larger organizations, with efforts typically centering large cities or state-wide and national projects (e.g. traffic stop analysis in Los Angeles, stop-and-frisk in New York City, and nationwide analyses of police use of force \cite{rossResolutionApparentParadoxes2018}, police shootings \cite{rossMultiLevelBayesianAnalysis2015} \cite{siegelRacialDisparitiesFatal}, and traffic stops \cite{tillyerImpactDriversRace2013}). However, according to the 2020 U.S. Census, about 76\% of the approximately 19,500 incorporated places in the U.S. had fewer than 5,000 people. Of those, almost 42\% had fewer than 500 people \cite{uscensusbureauLatestCityTown}. Smaller areas face barriers in accessing the data and support needed to analyze police data, which we seek to address in the SToPA Lab.

\subsection{Perspective}
\textbf{\textit{What is the context? Who in the community is working with and identifying the changes needed? Who in the community is working to make that change happen?}} \textbf{\textit{How did the knowledge of all research team members [including community members] guide the project in various directions?}}

The SToPA initiative began in Williamstown, MA. In 2020, a lawsuit was filed against the town of Williamstown because of alleged  racial discrimination and sexual abuse in the Williamstown Police Department \citep{cohen2020williamstown}. The Williamstown Racial Justice and Police Reform (WRJPR) group was formed in response. WRJPR included four working groups: social work, education, digital presence, and the data science group that eventually grew into SToPA \cite{williamstownracialjusticeandpolicereformWRJPRHomepage2021}.

\textbf{Place and context-based} Williamstown is a New England town that is home to Williams College, a private, selective liberal arts institution. It is home to approximately 7,500 residents, nearly 82\% of whom are white, with relatively high median income and low poverty rates. At the local level, Williamstown is sufficiently large to house its own service departments and open-meeting format government. The Williamstown Police Department (WPD) itself employs around 12 full-time officers on patrol at this time. Though there is a history of misconduct within the department, WPD came under severe scrutiny in August 2020 when a WPD officer filed a lawsuit against the town, its police chief, and town manager for racism and sexual harassment in the workplace \citep{cohen2020williamstown}. The resulting investigation unearthed a plethora of wrongdoing and abuses that were perpetrated for more than a decade in the department \citep{yang2022independent}. These abuses include evidence of Nazi idolatry among officers, illegal surveillance of community members, and evidence of rampant sexual and racist discrimination that was “initiated, participated in, and tolerated” within WPD \cite{joIndependent2022}. The lawsuit also shed light on the plaintiff, who had a history of misconduct and abuse dating as far back as the 1990s. This conduct within the WPD is the community context which motivated our work.

\subsection{Process} 
\textit{\textbf{How was the process of collaboration designed to enact the collaborative framework?}}

\textbf{Build mutual trust and respect, Relational and humanized} The SToPA project began with participatory data acquisition and cleaning by  a collective group of Williamstown residents, Williams College students, and off-site allied scholars. Here, we highlight the overlapping identities of researcher and community partner.  Multiple FOIA requests were filed which resulted in manually scanning thousands of printed pages of WPD data. Team members (consisting of community organizers, data scientists, statisticians, and sociological researchers) also delved into the racial and socioeconomic historical context of WPD's practices to gain greater insight into policing disparities. In one weekend, the data was cleaned to a minimally-usable state. Yet the newly cleaned data invited numerous questions regarding policing standards in Williamstown. After our initial data processing and analysis of data in Williamstown, we began to create a framework for community members from other localities to analyze their own data. This toolkit centers community-driven analysis by documenting and developing the tools that may be necessary to other communities interested in questions about policing. The toolkit also emphasizes the trust and relationships that must be developed through this analysis through a recommended team structure including people with a variety of skill sets and close knowledge of community context.

\subsection{Power}
\textit{\textbf{What aspects of the project were designed to provide communities with the power to take on this work? What additional capacity is enabled by this work?}}

\textbf{Mutually beneficial with coeducation} Members of SToPA later met with the WRJPR organizers to gain further context, listen to community experiences, and uncover which aspects of the data were most interesting to Williamstown residents. In this meeting, questions like “Is there a change in ticketing during the third shift?” and “Are there more traffic stops in the White Oaks [a predominantly Black] neighborhood?” were frequently posited. SToPA members soon realized that the addressing of such questions required experiential knowledge only community members could offer. The SToPA Lab was then formed to share insights on the community-identified issues and more broadly provide a reusable toolkit for other locales to conduct similar analyses using their own data.

\textbf{Democratic with shared ownership} The work of the SToPA Lab shows that it is possible to obtain, analyze, and interpret findings from policing data from small towns in order to enable data-driven advocacy. The software produced by the research lab is openly available for use, reuse and contributions on Github \cite{stopa_github}. However, SToPA Lab participants alone cannot analyze data from all of the small cities and towns in the country, and they do not have the proper community context to be able to do so. Though Williamstown was the initial focus of study, the work was undertaken with the intention of facilitating analysis of bias in policing in similar small towns. It was extended by residents of Rochester, NY, and Durham, NC, for use in their communities. At IMSI, we worked to develop a toolkit that would help scale these activities. Our research team included community activists and mathematicians,  including math researchers who have previously collected and cleaned data in collaboration with local community activist organizations.

In addition to SToPA’s origins in a Williamstown community activist group, SToPA team member Clarissa Aché Cabello analyzed policing in Durham, NC, and co-author KAİ has worked with community organization Decarcerate Utah to analyze Salt Lake City Police Department data. SToPA collaborators at IMSI also included local activists who were on the ground and most impacted by police bias. We worked with LTG, a queer and trans activist, researcher, and Harvard undergraduate student based in Nashville.

\begin{quote}
    \textbf{LTG}: I am a disabled, queer, nonbinary, trans, Black femme Pan-African LGBTQ+ community organizer. I am the Head of Justice with Queer Youth Assemble, a 2024 Freedom from Fear Fellow with Southerners on New Ground, a National Board Member of Queer Trans Project, and the National Community Demands Lead for the 2024 Queer Unity March. My research interests include intimate partner violence in Black queer relationships, the Black Queer Identity Matrix, and carceral abolitionist mental health care. I participated in IMSI the year before attending Harvard College as an undergraduate, months after graduating high school and relatively early in my organizing journey. My engagement in community organizations and initiatives greatly informed my understanding of community participatory action and its role in combating the criminal justice system. I have witnessed police violence in my communities and neighborhoods. And I have seen the power resourced community action has had in the struggle for change. Centering lived experiences and community needs in data science research is integral to aiding on-ground organizers in their fight for collective liberation.
\end{quote}

\begin{quote}
\textbf{Kİ}: I am a nonbinary, queer person of Middle Eastern descent who grew up witnessing the use of policing and state-sanctioned discrimination to target my community as “terrorists.” In the present, trans and nonbinary people are being criminalized for merely existing; all of my identities inspired me to participate in this work. My main contribution during my time at IMSI was sharing my experience working with raw data from my home of Salt Lake City’s Police Department (SLCPD) in an undergraduate research project, which required the use of RStudio to transform the raw data to remove duplicates, data cells listed as “Not Applicable”, typos, “Unknown” racial classifications made by police, conversion of dates of birth to ages, aligning racial categories used in the federal American Community Survey with those used by SLCPD, and more. Additionally, my previous experience organizing for environmental, reproductive, and LGBTQ+ justice, all of which emphasized the intersectional role of race, ethnicity, and other identity categories, helped me appreciate some of the decisions that local organizers may wish to make related to cleaning, transforming, and analyzing the data. In particular, my previous organizing experience emphasized the role of state-sanctioned oppression and required me to learn about the history and structural racism that can be embedded in state and local power structures. 

The SToPA toolkit is an abstraction of the work that is currently being conducted on a series of towns (Williamstown, Durham, Salt Lake City, Nashville and others) to be relevant and useful in conducting a more general research process in a variety of geographic settings. Through the toolkit, we developed generalizable methods to obtain the data, extract useful information from the original source, analyze it using statistical or visual methods, and use it to inform restorative justice efforts in a given community context.
\end{quote}

\begin{quote}
\textbf{AM}: By far the largest challenge we faced with the Williamstown data was getting it into a usable structure. Williamstown’s police records were stored in extremely dated technology, and the only format they could provide was thousands of printed pages, rather than any electronic format. There are so many small towns that likely face similar challenges– even if an activist or researcher successfully gets access to data, the delivered format of it may be entirely unusable. Much of the work we did in the first year was coding up a clean transformation from scanned PDF images to structured text data. Our hope is that this process provides a useful example to others, and documentation is crucial to its reproducibility.
\end{quote}

\begin{quote}
\textbf{SH}: My main contribution to the toolkit was providing documentation and examples of how to extract and transform data. Without the work that had been done with Williamstown and Durham I would simply not have been able to do this, as they provided concrete examples as to what data provided by police departments could look like. The code they produced was also important for my work, as I took their specific examples and worked to abstract it so that future groups would not have to start from scratch.
\end{quote}

\begin{quote}
\textbf{CK}: My previous experiences as a community organizer and as a statistician/data scientist made me conscious of many of the different use cases and data challenges that should be addressed by the toolkit. For example, policing data in Minneapolis is relatively clean and easy to access but is spatially privatized in ways that may not be obvious if a user is considering only individual events (and not visualizing the data spatially). However, this seemingly minor process of spatial privatization may have important implications for statistical analysis, for example when affiliating policing outcomes with census data. It is important to be aware of these seemingly minor modifications to the original data. Even when data is readily available, challenges arise because it often requires preprocessing before statistical analysis is appropriate, and the amount of publicly available data documentation varies. Incident-level data with spatial information is also often not available at all for many small towns. Through conversations with community members and my own organizing experiences, I tried to be mindful of the various use cases to make our resulting toolkit useful to as many communities as possible. I also created structured ways to shed light on some of the trickier aspects of publicly available data, that might have some complications that may not be immediately obvious.
\end{quote}

\subsection{Products and Relationships}
\textit{\textbf{How is the toolkit designed to guide users through a community-research partnership?}}

\begin{figure}[h]
    \centering
    \includegraphics[width = .9\textwidth]{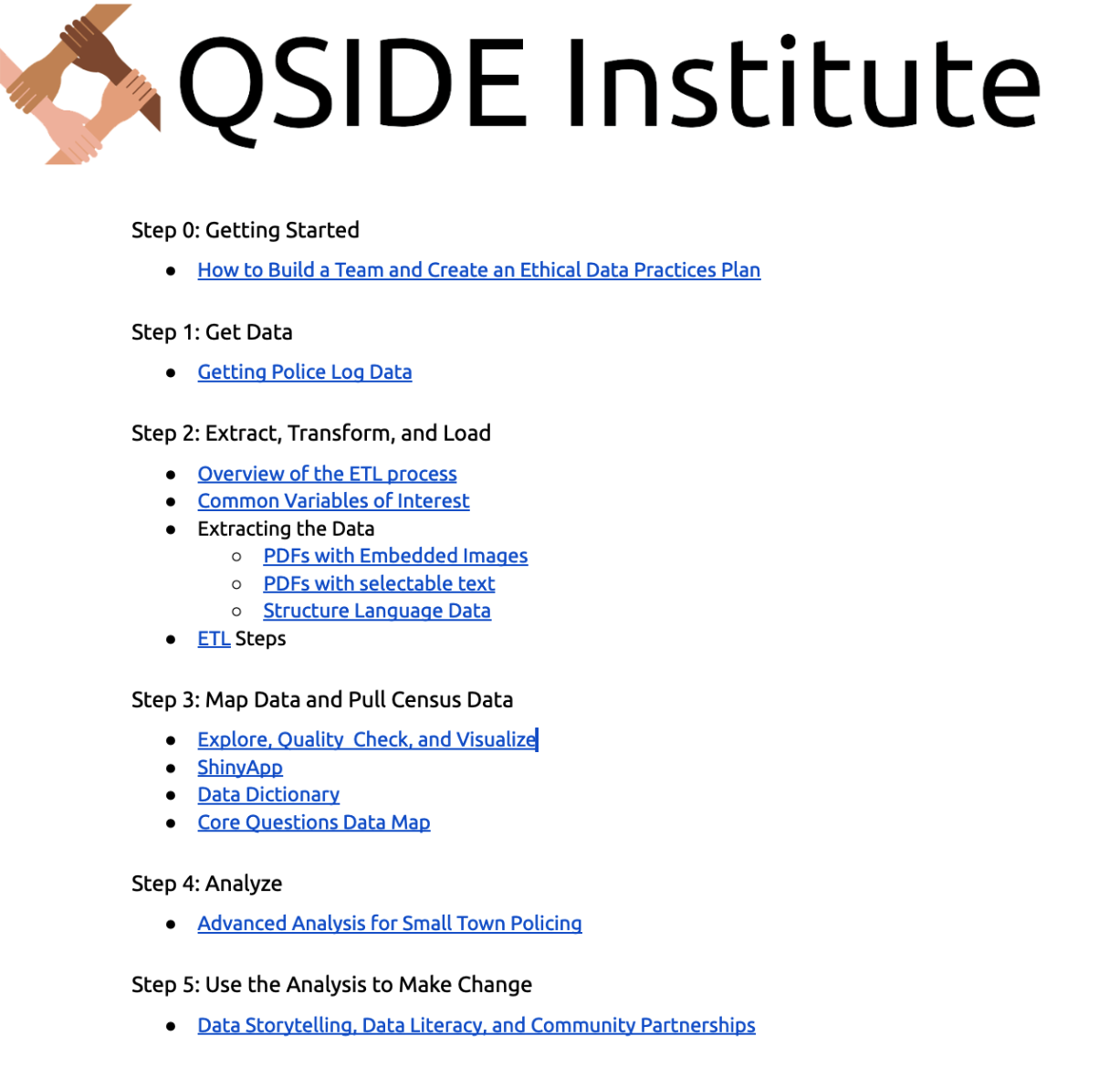}
    \caption{A \textit{draft} version of the interactive directory for the SToPA toolkit, where users will be able to find more information about how to obtain, clean, utilize, analyze, and act on findings from public data.}
    \label{fig:toolkit}
\end{figure}

\textbf{Developmental, sustained, and deep The draft SToPA toolkit exists as a document to help guide newcomers through the data collection and analysis process. Figure~\ref{fig:toolkit} shows a of the draft landing page of the toolkit}. We hope that the framework for community partnership is helpful to and utilized by future toolkit users. The SToPA toolkit guide highly recommends that, before starting collecting and cleaning data for analysis, researchers and community members join together in a community-based participatory action team (CBPAT) \cite{wilsonCommunityBasedParticipatoryAction2018}. For the purposes of the toolkit, CBPATs are considered to be groups of community members who possess a diversified set of intellectual and social skills. It is suggested that any CBPAT working to answer community-centered issues around policing have at least two types of participants: 
\begin{enumerate}
    \item \textbf{Community member(s) committed to police transparency and equity:} At least one member of the community – that is, someone who lives and/or works in the police jurisdiction– should be part of the team.
    \item \textbf{The Data Analyst:} This individual should have grounding in mathematics, statistics, and/or data science/analysis in order to accurately interpret quantitative findings. This will ensure that the toolkit can be adapted to the new setting.
\end{enumerate}

After these two participants are accounted for in a CBPAT, we also believe that additional team members could include:

\begin{enumerate}[resume]
    \item \textbf{The Social Changemaker:} This individual possesses skills in social organizing and community building. They likely will have extensive experience in community facilitation and outreach, communications skills, and intergroup mediation. Social Changemakers engage in work rooted in capturing and supporting community needs. This person would also directly engaging with other community members to ensure their voices are heard during the project. The CBPAT Social Changemaker is encouraged to accomplish this by working closely with other organizations in the community that are value-driven and aligned.
    \item \textbf{Social Scientist(s) and/or Attorney(s):} This individual should have a background in policing practice and research to ensure the findings are interpreted within appropriate legal and institutional frameworks.
    \item \textbf{Data Visualization Specialist(s)/Data Storyteller(s):} This individual would have an intermediate to advanced level of understanding in graphically representing relationships in data garnered from any collection efforts of the CBPAT. This Specialist’s expertise is especially useful in reframing visualizations for community storytelling and journalism.
\end{enumerate}

\textbf{Build capacity and agency} The toolkit was designed with all of these possible team members in mind, intending for teams to be  made up of largely local community members rather than outside researchers. Although some stages of the toolkit are more technical than others, each team member will be able to provide expertise throughout the process. We note that community context and social science expertise are particularly important when processing data, as documentation of policing data is often missing important details. These perspectives are also key in the analysis of the data and interpretation/dissemination of findings. The SToPA toolkit democratizes access to policing data analysis so that communities have the capacity to conduct their own analyses and the agency to drive local conversations about policing.  

\begin{quote} 
    \textbf{Kİ}: This toolkit was created in a scaffolded way so that groups without someone with data science knowledge could still obtain some results and that groups with someone who has some initial data science knowledge can go deeper into the analysis relevant to that community. Specifically, we aim to make the Small Town Policing Toolkit robust to many different use cases. For example, there may be towns or municipalities that have very structured data and can move quickly to the phase of analysis, while other cases will have extensive needs for data processing. There are also people or groups with different types of expertise that need to be equitably empowered through this process. For example, there may be folks with an expertise in community organizing that need more support in the data science side of the toolkit. There may also be experts in data science that need more support in the interpretation of the results in context. These ranging experiences, expertise sets, and data/situational contexts provide significant challenges for completing this kind of work. The toolkit is designed to be robust to many of these challenges.
\end{quote}

\section{Impact}
Both narratives from the communities that shaped the VECINA and SToPA projects share a sense of place, a sense of understanding needs and strengths from each other, and a desire to make a difference. But how do we know our work has been successful? A typical measure of impact in academic circles is a “product.” However, in this section, we intentionally broaden our discussion of impact beyond products and intertwine with personal narrative to illustrate these impacts.

VECINA in its fully realized form is a platform which will integrate traditional GIS tools with community knowledge, fuertes, and storytelling. However, the products at each stage in the collaboration were very different, and all important to the overall goal. SToPA’s work is to build a toolkit, but the research team was also building out long-term, sustainable plans for community engagement simultaneously. We emphasize that “Products” are not “Impact.” Products need to be relevant and accessible, but most importantly, they need people for real impact. Many of our team members have reflected on how this work has made an impact on us, on our work, and our hopes and dreams for the future.

\textbf{What impact has this made for you and/or your communities? What are your hopes and dreams for your own work moving forward?}

\begin{quote}
    \textbf{MJGP}: Finalmente, con la experiencia de tres años siendo parte de Nuevas Voces y el proyecto VECINA, con un destacable trabajo en equipo, dejo plasmada muchas ideas con información accesible, bilingüe y gratuita para la comunidad latianoamericana en general que vino a buscar mejores días y deseo cubrir sus necesidades hablando el mismo idioma¨ en un país extranjero lleno de oportunidades de progreso¨.
    \textbf{MJGP} (translation by JM): Finally, with the experience of three years being part of Nuevas Voces and the VECINA project, with remarkable teamwork, I have captured many ideas with accessible, bilingual and free information for the Latin American community in general who came looking for better days and wish meet their needs speaking the same language ``in a foreign country full of opportunities for advancement''.
\end{quote}

\begin{quote}
\textbf{RG}: I am still trying to wrap my head around all the ways I continue to be personally impacted, and all the ways WRWC and Nuevas Voces are expanding into new possibilities out of our participation on this project. As a long-time non-profit community builder with training in Critical Theory and International Human Rights and decades of boots on the ground work in the arts and environmental sectors, I have both felt and seen a lot of hopes and dreams, projects and partnerships, and collective capacity come and go, but I think something real, authentic and new is happening here. The model is strong, fostering shared empowerment, intentionality, commitment and sense of purpose, and the people are even stronger. 

I have faith in our collective ability to launch this site, situate it in the hands of the people in Olneyville, and secure continuous resources to expand its impact and reach. Hope is harder. I hope to see equity, justice and all aspects of ``Buen Vivir" thrive in Olneyville and nearby neighborhoods. I hope I continue to become a better and better ally as I learn from the experts - those with firsthand experience and those who have been raised in more awakened times. My dream is that through coming together we can overturn some of the ravages of racism and plunder one community, one project, one website at a time, and that in coming together we can share the biggest possible dream and make it happen. Time will tell.
\end{quote}

\begin{quote}
\textbf{CDE}: I went into the VECINA project wanting to give back to an important community in my life. I wanted to enact a meaningful process for all involved that was true to our philosophy of how communities should be engaged. I met and worked with wonderful people at WRWC and Nuevas Voces. And a milestone - I facilitated my first research group in Spanish! But what I did not expect was to see my colleagues at WRWC and Nuevas Voces happy and excited about data dashboards, working with mathematicians, and the future possibilities ahead. I believe this work in partnership with communities is also changing the perception about math and mathematicians in addition to making a difference with data science.
\end{quote}

\begin{quote}
\textbf{JEH}: The project has given me the opportunity to work with amazing mathematicians and amazing people in the community of Providence.  In my courses that I teach I have spent so much time talking about working with communities and now I have some evidence about the best practices.  I am so grateful to learn from our partners at WRWC and Nuevas Voces.  I am nervous about the future of the project, but I believe that VECINA is a sustainable project for the community.
\end{quote}

\begin{quote}
\textbf{SH}: In my time as a mathematical communicator and librarian I have long been focused on the human stories of mathematics, but most of the time this meant the humans who are doing the mathematics, not the ones who are being impacted by them. It is only in the last couple of years that I began to comprehend just how bijective a function exists between mathematics and humanity. My time at IMSI, both my closer working with the SToPA toolkit and interacting with the VECINA team, helped me really concretize this understanding, mathematics is both done by humans, has a tremendous impact on humans, and it is up to humans to make sure that that impact is positive. This means that it is up to all of us who do mathematics to truly think about work that we are doing beyond its beauty and elegance and abstract nature, reaching to consider what it could mean for the world beyond its ability to help lay the foundation for yet another proof. 
\end{quote}

\begin{quote}
    \textbf{VP}: I was a lawyer before becoming a mathematician, and I am always looking to combine these two experiences in a way to be of service. I could not imagine a more meaningful way to do this than working with community partners and collaborating with them to solve problems. In the work we did on the Providence public schools, we both read and interpreted statutes and conducted data-based inquiry. In the end, I hope we have a tool that residents of Olneyville find useful.
\end{quote}

\begin{quote}
\textbf{CK}: As a statistician, I am interested in making sure that questions from the community are driving the statistical analysis, rather than developing statistical methods in search of a dataset or problem that is relevant to it. This experience at IMSI and with the SToPA team has been a transformative experience, allowing me to work with folks that are not only motivated by real-world problems, but are involved in their communities and seeking to create positive change by co-creating that change with community members. I am excited to see how the work will be used with communities to translate data or ideas into insight and action.
\end{quote}

\begin{quote}
\textbf{AM}: While the work the SToPA research lab has accomplished so far is impactful, I am eager to connect more closely with the Williamstown community to transfer what we have and reconnect with our broader goals. Throughout this work, Dr. Bilal Ansari, an influential Williamstown community member, has provided invaluable feedback and guidance. He has expressed excitement for our results and findings, however, much of the work has been somewhat siloed because of its tedious, technical nature. Moving forward, my goal for this particular work is to finalize it in a way that the Williamstown community can explore and take ownership over the project. My vision for the future of the lab is to create a repository that is accessible to and used by a diversity of stakeholders: community members, data scientists and researchers alike. Aside from the community-facing impact, this work has provided me with a unique opportunity to advance my research program with incredible collaborators and apply my skills in an area that matches with my values. I find immense fulfillment in being able to connect with communities and share my knowledge to solve problems that are meaningful to people and their daily lives, especially ones that uplift the most vulnerable. 
\end{quote}

\begin{quote}
\textbf{Kİ}: As a mathematician trained in knot theory, my data science skills have been mostly self-taught, so the opportunity to work with experts across fields–including social scientists trained in analyzing the social structures that affect how our mathematics is perceived and used in the world, data scientists with experience in coding and choosing the best analysis techniques to tell the story behind the data, and perhaps most importantly, community members and activists most impacted by racialized policing–was illuminating. 

I believe it is important for academics who study fields with real-world impact–arguably all of us–to realize that academic expertise is only one form of knowledge, and that those most impacted by the phenomena we study have personal experience–a holistic form of knowledge and a perspective that is arguably broader than our own. I hope to broaden my own collaboration with community members in future work, realizing that academic training has epistemically devalued personal experience and combating my own internalized “ivory tower” ideologies.

My hope is that mathematics, as a field, is able to more greatly consider the impacts of its work on communities, particularly members of minoritized groups, and collaborate with those most impacted. Some mathematicians have historically distanced themselves from their research and the ways it is used–for example, NSA-funded cryptographic research and the development of algorithms to determine recidivism risk \cite{oneilWeaponsMathDestruction2016a}–and I hope that by collaborating with affected community members, mathematics can become a force for a fairer and more just world.
\end{quote}

\begin{quote}
\textbf{LTG}:   I contributed to the project by attempting to create avenues of communication between local activist groups and data science researchers as an on-ground organizer. However, because of my limited experience in organizing at the time of my participation in this program, the extent of this collaboration between researchers and activist groups on ground nationally within the SToPA IMSI Project was also limited. Instead, I relied on the direct and indirect collective guidance of scholar-activists and organizers from my community to inform how I engaged in the space.

 I hope our initial investment in the budding collaboration among Community Members, Data Analysts, and Social Changemakers within SToPA starts a movement of intentional and direct involvement of experienced, educated organizers and activists that center minoritized experiences. I dream of a network of civically engaged individuals of all education levels who find inspiration from our frameworks and toolsets to build larger initiatives where community members visualize and actualize their freedom dreams of safety and security free from policing.  As a young person and undergraduate student, I hope that the invitation of empowered youth to rigorous research and collaborative spaces in data science and participatory action continues beyond me. 
 \end{quote}

\section{Challenges and Opportunities}

The two projects presented here brought opportunities and challenges. We focus first on challenges that were shared across both projects. Even though there were different levels of community involvement at IMSI, both projects center community needs and kept this audience in mind as tools were developed. Community partners were involved in the research, design, and deployment of all parts of the projects. Securing consistent and ongoing funding for community partners and other research needs has also been a challenge for both teams.

\begin{quote}
\textbf{BAT}: VECINA has partnered with math institutes several times now, as well as Bates College and the RIOS institute. Working with the “ivory tower” has given us opportunities to work together, but also has introduced challenges on who is seen as a researcher. There was not always office space provided for community members and finding funding for travel and lodging, per diem, or stipends has not always been straightforward. It has been important to the team to not center every part of a sprint entirely at the institutes by making sure to travel to the community events to learn from them and be a part of the Nuevas Voces classes. As the group has grown, keeping communication channels open and inviting new community members into VECINA has been a priority as not everyone uses the same tools (email, text, etc) and we keep planning for scalability in the future.
\end{quote}

\begin{quote}
\textbf{JEH}: When working with the institutes, it has been a challenge to convince them and leadership to give resources to have community members be a part of the research team. Another challenge is convincing the institutes to be willing to open the doors to the community and recognize that community members are part of the mathematical community even if they don’t have a traditional mathematical degree. We’ll continue to work with them to help change perceptions. 
\end{quote}

We also note challenges and lessons we saw throughout the process of drafting the SToPA toolkit and during followup work being conducted now in the DSPACE research group. Developing the SToPA toolkit requires substantial time investment of team members, which creates challenges around maintaining relevance and sustaining community partnerships throughout the production process. For our junior faculty members, we also recognized that toolkit development may not align with traditional tenure and promotion criteria, which motivated the creation of our research wing in DSPACE to ensure the scholarly contributions are in formats beneficial to both community members and faculty. Additionally, we remain mindful that our current community partners in DSPACE are volunteers with limited capacity. We strive to engage them meaningfully and frequently while respecting the constraints on their time and resources.

There is great potential for researchers in mathematics to conduct meaningful and impactful analysis alongside community members. We believe that our framework and the three principles of power, perspective, and process can provide structure to support intentional and ethical collaboration. We offer descriptions of our experiences and individual reflections to encourage and equip other researchers to pursue this vital work.

\bibliography{references}

\vfill\eject

\end{document}